\documentclass[10pt]{article}
\usepackage{cite}
\usepackage{mathrsfs}
\usepackage{amsfonts}
\usepackage{amsmath}
\usepackage{amsfonts,amssymb}
\usepackage{dsfont}
\usepackage{curves}
\usepackage{mathrsfs}
\usepackage{pifont}
\usepackage{amssymb}
\allowdisplaybreaks

\numberwithin{equation}{section}

\date{}

\def\BigRoman{\uppercase\expandafter{\romannumeral\number\count 255 }}
\def\Romannumeral{\afterassignment\BigRoman\count255=}

\begin{document}
\title{Distance spectral radius for a graph to be $k$-critical with respect to $[1,b]$-odd factor
%\thanks{Supported by }
}
\author{\small Sufang Wang$^{1}$, Wei Zhang$^{2}$\footnote{Corresponding
author. E-mail address: wangsufangjust@163.com (S. Wang), zw\_wzu@163.com (W. Zhang).}\\
\small $1$.  School of Public Management, Jiangsu University of Science and Technology,\\
\small Zhenjiang, Jiangsu 212100, China\\
\small $2$.  School of Economics and Management, Wenzhou University of Technology,\\
\small Wenzhou, Zhejiang 325000, China
}

\maketitle
\begin{abstract}
\noindent Let $G$ be a connected graph, and let $b$ and $k$ be two positive integers with $b\equiv1$ (mod 2). A $[1,b]$-odd factor of $G$ is
a spanning subgraph $F$ of $G$ with $d_F(v)\equiv1$ (mod 2) and $1\leq d_F(v)\leq b$ for every $v\in V(G)$. A graph $G$ is called $k$-critical
with respect to $[1,b]$-odd factor if $G-X$ contains a $[1,b]$-odd factor for every $X\subseteq V(G)$ with $|X|=k$. Let $\mathcal{D}(G)$ denote
the distance matrix of $G$. The largest eigenvalue of $\mathcal{D}(G)$, denoted by $\mu(G)$, is called the distance spectral radius of $G$. In
this paper, we prove an upper bound for $\mu(G)$ in a $(k+1)$-connected graph $G$ which guarantees $G$ to be $k$-critical with respect to $[1,b]$-odd
factor.
\\
\begin{flushleft}

(2020) Mathematics Subject Classification: 05C50, 05C70, 90B99
\end{flushleft}
\end{abstract}

\section{Introduction}

Throughout this paper, we consider finite and undirected graphs without multiple edges or loops. Let $G$ be a graph with vertex set
$V(G)=\{v_1,v_2,\ldots,v_n\}$ and edge set $E(G)$. The order of $G$ is the number $n=|V(G)|$ of its vertices. A graph $G$ is called trivial
if its order is 1. The degree of a vertex $v$ in $G$, denoted by $d_G(v)$, is the number of vertices adjacent to $v$. Let $o(G)$ denote the
number of odd components in $G$. For any $S\subseteq V(G)$, we denote by $G[S]$ and $G-S$ the subgraphs of $G$ induced by $S$ and $V(G)\setminus S$,
respectively. Let $K_n$ denote the complete graph of order $n$.

For $v_i,v_j\in V(G)$, the distance between $v_i$ and $v_j$, denoted by $d_{ij}$, is the length of a shortest path from $v_i$ to $v_j$. The
distance matrix of $G$ is defined as $\mathcal{D}(G)=(d_{ij})_{n\times n}$. The largest eigenvalue of distance matrix $\mathcal{D}(G)$, which
is also called the spectral radius of $\mathcal{D}(G)$, is denoted by $\mu(G)$. By virtue of the Perron-Frobenius theorem, $\mu(G)$ is always
positive (unless $G$ is trivial), and we call $\mu(G)$ the distance spectral radius in $G$. Graham and Pollack \cite{GP} described a connection
between the number of negative distance eigenvalues and the addressing problem in data communication system, which is regarded as the original
study of distance eigenvalues. Since then, the study of distance eigenvalues of graphs has become a highly concerned research topic, and this
topic has attracted more and more attention. Aouchiche and Hansen \cite{AH} verified lower and upper bounds on the distance spectral radius of
a connected graph using proximity and remoteness. Mojallal and Hansen \cite{MH} established a relation between proximity and the third largest
distance eigenvalue in a connected graph. Zhang and Lin \cite{ZL} characterized all trees with three distinct distance eigenvalues. Zhou and Wu
\cite{ZW} proved an upper bound on the distance spectral radius for a connected graph to contain a spanning $k$-tree. Zhou \cite{Zt} provided
an upper bound on the distance spectral radius of a connected graph $G$ to guarantee that $G$ is $\frac{1}{t}$-tough and $t$-tough, respectively.

Let $G_1$ and $G_2$ be two vertex-disjoint graphs. The disjoint union of $G_1$ and $G_2$ is denoted by $G_1\cup G_2$. The join of $G_1$ and $G_2$
is denoted by $G_1\vee G_2$, which is the graph obtained from $G_1\cup G_2$ by adding all possible edges between $V(G_1)$ and $V(G_2)$.

A $[1,b]$-odd factor of $G$ is a spanning subgraph $F$ of $G$ with $d_F(v)\equiv1$ (mod 2) and $1\leq d_F(v)\leq b$ for every $v\in V(G)$, where
$b$ is an odd positive integer. A $[1,1]$-odd factor is also called a perfect matching. A graph $G$ is called $k$-critical with respect to
$[1,b]$-odd factor if $G-X$ contains a $[1,b]$-odd factor for every $X\subseteq V(G)$ with $|X|=k$. A $k$-critical graph with respect to $[1,1]$-odd
factor is simply called a $k$-factor-critical graph.

Tutte \cite{Tutte} proposed a criterion for a graph with a perfect matching. Enomoto \cite{E} gave a toughness condition for the existence of a
perfect matching in a graph. Anderson \cite{Anderson} established a connection between binding number and a perfect matching in a graph. Brouwer
and Haemers \cite{BH} provided an eigenvalue condition for the existence of a perfect matching in a regular graph. O \cite{Os} characterized a
perfect matching of a graph with respect to its size and spectral radius. Zhang and Lin \cite{ZL2} proposed a distance spectral radius condition
which ensures a graph to contain a perfect matching. Zhou and Zhang \cite{ZZ} established a relationship between signless Laplacian spectral radius
and a perfect matching in a graph. Zhou, Xu and Sun \cite{ZXS} obtained some results on $[1,b]$-factors in graphs. Amahashi \cite{A} presented a
characterization for a graph with a $[1,b]$-odd factor. Cui and Kano \cite{CK} provided a neighborhood condition for the existence of a $[1,b]$-odd
factor in a graph. Egawa, Kano and Yan \cite{EKY} investigated $[1,b]$-odd factors in graphs and obtained some results on the existence of
$[1,b]$-odd factors in graphs. Lu, Wu and Yang \cite{LWY} got some sufficient conditions by virtue of eigenvalues for the existence of $[1,b]$-odd
factors. Kim, O, Park and Ree \cite{KOPR} posed an eigenvalue condition for a graph to have a $[1,b]$-odd factor. Zhou and Liu \cite{ZL3} showed two
sufficient conditions for a connected graph to contain a $[1,b]$-odd factor via its size and spectral radius. Favaron \cite{Favaron} proved a
necessary and sufficient condition for a graph to be $k$-factor-critical. Zhou, Sun and Zhang \cite{ZSZ} provided a lower bound for the spectral
radius in a graph which guarantees the graph to be $k$-factor-critical. Kano and Matsuda \cite{KM} established a criterion for a graph to be
$k$-critical with respect to $[1,b]$-odd factor. For more results on factor critical graphs, we refer the reader to \cite{ACA,PS,YH,ZZL,ZPX1,Zr}.
Lots of researchers \cite{Ws,Wu,WZ,ZZL1,Zs,Zs1} presented some sufficient spectral conditions for the existence of graph factors.

Motivated by \cite{KM,ZL2} directly, it is natural and interesting to propose a distance spectral radius condition which guarantees a graph to be
$k$-critical with respect to $[1,b]$-odd factor. In this paper, we investigate the connection between distance spectral radius and a $k$-critical
graph with respect to $[1,b]$-odd factor, and give a distance spectral radius condition for a connected graph to be $k$-critical with respect to
$[1,b]$-odd factor. Our main result is shown as follows.

\medskip

\noindent{\textbf{Theorem 1.1.}} Let $b,k$ and $n$ be three positive integers with $b\equiv1$ (mod 2) and $n\equiv k$ (mod 2), and let $G$ be a
$(k+1)$-connected graph of order $n\geq\frac{b^{2}+2bk+5b+2k+4}{b}$. If $G$ satisfies
$$
\mu(G)\leq\mu(K_{k+1}\vee(K_{n-k-b-2}\cup(b+1)K_1)),
$$
then $G$ is $k$-critical with respect to $[1,b]$-odd factor, unless $G=K_{k+1}\vee(K_{n-k-b-2}\cup(b+1)K_1)$.

\medskip

\section{Some preliminaries}

Let $f:V(G)\rightarrow\{1,3,5,\ldots\}$ be a function. A $(1,f)$-odd factor of a graph $G$ is a spanning subgraph $F$ of $G$ with $d_F(v)\equiv1$
(mod 2) and $1\leq d_F(v)\leq f(v)$ for every $v\in V(G)$. Let $k$ be a positive integer. A graph $G$ of order $n\geq k+2$ is $k$-critical with
respect to $(1,f)$-odd factor if $G-X$ has a $(1,f)$-odd factor for every $X\subseteq V(G)$ with $|X|=k$. Kano and Matsuda \cite{KM} provided a
characterization for a graph to be $k$-critical with respect to $(1,f)$-odd factor.

\medskip

\noindent{\textbf{Lemma 2.1}} (Kano and Matsuda \cite{KM}). Let $k\geq1$ be an integer, and let $G$ be a graph of order $n\geq k+2$. Then $G$ is
$k$-critical with respect to $(1,f)$-odd factor if and only if for any $S\subseteq V(G)$ with $|S|\geq k$,
$$
o(G-S)\leq\sum\limits_{v\in S}{f(v)}-\max\Big\{\sum\limits_{v\in X}{f(v)}: X\subseteq S, |X|=k\Big\}.
$$

\medskip

\noindent{\textbf{Lemma 2.2}} (Minc \cite{Mn}). Let $G$ be a connected graph with two nonadjacent vertices $u,v\in V(G)$. Then
$$
\mu(G)>\mu(G+uv),
$$
where $G+uv$ denotes the graph obtained from $G$ by adding an edge to connect $u$ and $v$.

\medskip

The Wiener index of a connected graph $G$ is defined by $W(G)=\sum\limits_{i<j}d_{ij}$. The following lemma can be easily derived by the Rayleigh
quotient \cite{HJ}.

\medskip

\noindent{\textbf{Lemma 2.3.}} Let $G$ be a connected graph of order $n$. Then
$$
\mu(G)=\max\limits_{X\neq\textbf{0}}\frac{X^{T}\mathcal{D}(G)X}{X^{T}X}\geq\frac{\textbf{1}^{T}\mathcal{D}(G)\textbf{1}}{\textbf{1}^{T}\textbf{1}}=\frac{2W(G)}{n},
$$
where $\textbf{1}=(1,1,\ldots,1)^{T}$.

\medskip

\noindent{\textbf{Lemma 2.4}} (Zheng, Li, Luo and Wang \cite{ZLLW}). Let $n=\sum\limits_{i=1}^{t}n_i+s$. If $n_1\geq n_2\geq\cdots\geq n_t\geq p\geq1$
and $n_1<n-s-p(t-1)$, then
$$
\mu(K_s\vee(K_{n_1}\cup K_{n_2}\cup\cdots\cup K_{n_t}))>\mu(K_s\vee(K_{n-s-p(t-1)}\cup(t-1)K_p)).
$$

\medskip

Let $M$ be an $n\times n$ real matrix and $V=\{1,2,\ldots,n\}$. Given a partition $\pi: V=V_1\cup V_2\cup\cdots\cup V_r$, the matrix $M$ can be
correspondingly partitioned as
\begin{align*}
M=\left(
  \begin{array}{cccc}
    M_{11} & M_{12} & \cdots & M_{1r}\\
    M_{21} & M_{22} & \cdots & M_{2r}\\
    \vdots & \vdots & \ddots & \vdots\\
    M_{r1} & M_{r2} & \cdots & M_{rr}\\
  \end{array}
\right),
\end{align*}
where $M_{ij}$ denotes the submatrix (block) of $M$ formed by rows in $V_i$ and columns in $V_j$. Let $b_{ij}$ be the average row sum of $M_{ij}$.
The quotient matrix of $M$ with respect to $\pi$ is defined as the matrix $B_{\pi}=(b_{ij})_{r\times r}$. The partition $\pi$ is called equitable
if every block $M_{ij}$ of $M$ admits constant row sum $b_{ij}$.

\medskip

\noindent{\textbf{Lemma 2.5}} (You, Yang, So and Xi \cite{YYSX}). Let $M$ be a real symmetric matrix with an equitable partition $\pi$, and let
$B_{\pi}$ be the corresponding quotient matrix. Then the eigenvalues of $B_{\pi}$ are also eigenvalues of $M$. Furthermore, if $M$ is nonnegative
and irreducible, then the largest eigenvalues of $M$ and $B_{\pi}$ are equal.

\section{The proof of Theorem 1.1}

\noindent{\it Proof of Theorem 1.1.} Suppose to the contrary that a $(k+1)$-connected graph $G$ is not $k$-critical with respect to $[1,b]$-odd
factor. By virtue of Lemma 2.1, we possess $o(G-S)>b|S|-bk$ for some nonempty subset $S\subseteq V(G)$ with $|S|\geq k$. If $|S|=k$, then
$0=o(G-S)>0$, a contradiction. Therefore, $|S|\geq k+1$. Write $|S|=s\geq k+1$ and $o(G-S)=q$. Since $n\equiv k$ (mod 2), $q$ and $s-k$ admit the
same parity. Together with $b\equiv1$ (mod 2), $q$ and $bs-bk$ possess the same parity. Hence, we conclude $q\geq bs-bk+2$. Then $G$ is a spanning
subgraph of $G_1=K_s\vee(K_{n_1}\cup K_{n_2}\cup\cdots\cup K_{n_{bs-bk+2}})$ for some positive odd integers $n_1\geq n_2\geq\cdots\geq n_{bs-bk+2}$
with $\sum\limits_{i=1}^{bs-bk+2}n_i=n-s$. By means of Lemma 2.2, we get
\begin{align}\label{eq:3.1}
\mu(G)\geq\mu(G_1),
\end{align}
where the equality holds if and only if $G=G_1$. Set $G_2=K_s\vee(K_{n-(b+1)s+bk-1}\cup(bs-bk+1)K_1)$, where $n\geq(b+1)s-bk+2$. By Lemma 2.4, we
admit
\begin{align}\label{eq:3.2}
\mu(G_1)\geq\mu(G_2),
\end{align}
with equality holding if and only if $(n_1,n_2,\ldots,n_{bs-bk+2})=(n-(b+1)s+bk-1,1,\ldots,1)$. The quotient matrix of $\mathcal{D}(G_2)$ with the
partition $V(G_2)=V(K_s)\cup V(K_{n-(b+1)s+bk-1})\cup V((bs-bk+1)K_1)$ is equal to
\begin{align*}
B=\left(
  \begin{array}{ccc}
    s-1 & n-(b+1)s+bk-1 & bs-bk+1\\
    s & n-(b+1)s+bk-2 & 2bs-2bk+2\\
    s & 2n-2(b+1)s+2bk-2 & 2bs-2bk\\
  \end{array}
\right).
\end{align*}
Then the characteristic polynomial of $B$ is
\begin{align*}
f_B(x)=&x^{3}-(n+bs-bk-3)x^{2}\\
&+((2bk-2bs-5)n+(2b^{2}+3b)s^{2}-(4b^{2}k+3bk-3b-3)s+2b^{2}k^{2}-3bk+6)x\\
&-(b^{2}+b)s^{3}+(bn+2b^{2}k+bk+2b^{2}+b-1)s^{2}\\
&+((1-2b-bk)n-b^{2}k^{2}-4b^{2}k-bk+4b+2)s\\
&+(2bk-4)n+2b^{2}k^{2}-4bk+4.
\end{align*}
Note that the partition $V(G_2)=V(K_s)\cup V(K_{n-(b+1)s+bk-1})\cup V((bs-bk+1)K_1)$ is equitable. Then it follows from Lemma 2.5 that $\mu(G_2)$
is the largest root of $f_B(x)=0$.

Let $G_*=K_{k+1}\vee(K_{n-k-b-2}\cup(b+1)K_1)$. The quotient matrix of $\mathcal{D}(G_*)$ in terms of the partition $V(G_*)=V(K_{k+1})\cup V(K_{n-k-b-2})\cup V((b+1)K_1)$
can be expressed as
\begin{align*}
B_*=\left(
  \begin{array}{ccc}
    k & n-k-b-2 & b+1\\
    k+1 & n-k-b-3 & 2b+2\\
    k+1 & 2n-2k-2b-4 & 2b\\
  \end{array}
\right).
\end{align*}
By a direct computation, the characteristic polynomial of $B_*$ is
\begin{align*}
f_{B_*}(x)=&x^{3}-(n+b-3)x^{2}+(-(2b+5)n+2b^{2}+3bk+6b+3k+9)x\\
&+(bk+k-b-3)n-b^{2}k+b^{2}-bk^{2}-bk+4b-k^{2}+5.
\end{align*}
Since the partition $V(G_*)=V(K_{k+1})\cup V(K_{n-k-b-2})\cup V((b+1)K_1)$ is equitable, from Lemma 2.5, the largest root, say $\theta$, of
$f_{B_*}(x)=0$ is equal to $\mu(G_*)$. That is to say, $\mu(G_*)=\theta$ and $f_{B_*}(\theta)=0$.

Recall that $G_2=K_s\vee(K_{n-(b+1)s+bk-1}\cup(bs-bk+1)K_1)$. If $s=k+1$, then $G_2=K_{k+1}\vee(K_{n-k-b-2}\cup(b+1)K_1)$. Combining this with
(\ref{eq:3.1}) and (\ref{eq:3.2}), we obtain
$$
\mu(G)\geq\mu(K_{k+1}\vee(K_{n-k-b-2}\cup(b+1)K_1)),
$$
with equality if and only if $G=K_{k+1}\vee(K_{n-k-b-2}\cup(b+1)K_1)$. Observe that $G=K_{k+1}\vee(K_{n-k-b-2}\cup(b+1)K_1)$ is not $k$-critical
with respect to $[1,b]$-odd factor, a contradiction. Next, we consider $s\geq k+2$.

By plugging the value $\theta$ into $x$ of $f_B(x)-f_{B_*}(x)$, we get
\begin{align}\label{eq:3.3}
f_B(\theta)=f_B(\theta)-f_{B_*}(\theta)=(s-k-1)g(\theta),
\end{align}
where $g(\theta)=-b\theta^{2}+(-2bn+(2b^{2}+3b)s-2b^{2}k+2b^{2}+6b+3)\theta-(b^{2}+b)s^{2}+(bn+b^{2}k+b^{2}-1)s+(1-b)n-2b^{2}k-bk+b^{2}+4b-k+1$.

According to Lemma 2.3 and $n\geq\frac{b^{2}+2bk+5b+2k+4}{b}$, we conclude
\begin{align}\label{eq:3.4}
\theta=&\mu(G_*)\geq\frac{2W(G_*)}{n}\nonumber\\
=&\frac{n^{2}+(2b+1)n-b^{2}-2bk-5b-2k-4}{n}\nonumber\\
\geq&n+b+1.
\end{align}

The symmetry axis of $g(\theta)$ is $\theta=\frac{-2bn+(2b^{2}+3b)s-2b^{2}k+2b^{2}+6b+3}{2b}$, which implies that $g(\theta)$ is strictly
decreasing in the interval $[\frac{-2bn+(2b^{2}+3b)s-2b^{2}k+2b^{2}+6b+3}{2b},+\infty)$. In view of (\ref{eq:3.4}), $s\geq k+2$ and
$n\geq(b+1)s-bk+2$, we have
$$
\frac{-2bn+(2b^{2}+3b)s-2b^{2}k+2b^{2}+6b+3}{2b}<n+b+1\leq\theta.
$$
Thus, we deduce
\begin{align*}
g(\theta)\leq&g(n+b+1)\\
=&-b(n+b+1)^{2}+(-2bn+(2b^{2}+3b)s-2b^{2}k+2b^{2}+6b+3)(n+b+1)\\
&-(b^{2}+b)s^{2}+(bn+b^{2}k+b^{2}-1)s+(1-b)n-2b^{2}k-bk+b^{2}+4b-k+1\\
=&-(b^{2}+b)s^{2}+((2b^{2}+4b)n+2b^{3}+b^{2}k+6b^{2}+3b-1)s-3bn^{2}\\
&+(-2b^{2}k-2b^{2}+b+4)n-2b^{3}k+b^{3}-4b^{2}k+7b^{2}-bk+12b-k+4\\
\leq&-(b^{2}+b)\Big(\frac{n+bk-2}{b+1}\Big)^{2}+((2b^{2}+4b)n+2b^{3}+b^{2}k+6b^{2}+3b-1)\Big(\frac{n+bk-2}{b+1}\Big)\\
&-3bn^{2}+(-2b^{2}k-2b^{2}+b+4)n-2b^{3}k+b^{3}-4b^{2}k+7b^{2}-bk\\
&+12b-k+4 \ \ \ \ \ \ \Big(\mbox{since} \ s\leq\frac{n+bk-2}{b+1}\Big)\\
=&\frac{1}{b+1}(-b^{2}n^{2}+(b^{2}k+b^{2}+4b+3)n+b^{4}+4b^{3}+7b^{2}-3bk+6b-k+6)\\
\leq&\frac{1}{b+1}\Big(-b^{2}\Big(\frac{b^{2}+2bk+5b+2k+4}{b}\Big)^{2}+(b^{2}k+b^{2}+4b+3)\Big(\frac{b^{2}+2bk+5b+2k+4}{b}\Big)\Big)\\
&+\frac{b^{4}+4b^{3}+7b^{2}-3bk+6b-k+6}{b+1} \ \ \ \ \ \ \Big(\mbox{since} \ n\geq \frac{b^{2}+2bk+5b+2k+4}{b}\Big)\\
=&\frac{1}{b(b+1)}(-3b^{4}k-5b^{4}-2b^{3}k^{2}-17b^{3}k-17b^{3}-6b^{2}k^{2}-25b^{2}k\\
&-7b^{2}-4bk^{2}-3bk+21b+6k+12)\\
<&0.
\end{align*}
Combining this with (\ref{eq:3.3}) and $s\geq k+2$, we obtain
$$
f_B(\theta)=(s-k-1)g(\theta)<0,
$$
which leads to
\begin{align}\label{eq:3.5}
\mu(G_2)>\theta=\mu(G_*)=\mu(K_{k+1}\vee(K_{n-k-b-2}\cup(b+1)K_1)).
\end{align}
In terms of (\ref{eq:3.1}), (\ref{eq:3.2}) and (\ref{eq:3.5}), we deduce
$$
\mu(G)\geq\mu(G_1)\geq\mu(G_2)>\mu(K_{k+1}\vee(K_{n-k-b-2}\cup(b+1)K_1)),
$$
which is a contradiction to $\mu(G)\leq\mu(K_{k+1}\vee(K_{n-k-b-2}\cup(b+1)K_1))$. This completes the proof of Theorem 1.1. \hfill $\Box$

\medskip

\section*{Declaration of competing interest}

\medskip

The authors declares that they have no known competing financial interests or personal relationships that could have appeared to influence the work
reported in this paper.

\section*{Data availability}

\medskip

No data was used for the research described in the article.

\medskip

%\section*{Acknowledgments}

\end{document}